\documentclass{article}
\usepackage{epsfig}
\usepackage{latexsym}
\usepackage{amsfonts}
\def\card{\mathrm{card}\,}
\def\Rea{\mathrm{Re}\,}
\def\Ima{\mathrm{Im}\,}
\def\C{\mathbf{C}}
\def\id{\mathrm{id}}
\def\N{\mathfrak{N}}
\def\Fr{\mathfrak{Fr}}
\begin{document}
\title{Densities in Fabry's theorem}
\author{Alexandre Eremenko\thanks{Supported
by NSF grant DMS-0555279}}
\maketitle
\begin{abstract}
Fabry's theorem on the singularities of power series 
is improved: the maximum
density in the assumptions of this theorem is
replaced by an 
interior density of Beurling--Malliavin type.

MSC classes: 30B10, 30B40.
\end{abstract}

\noindent
{\bf 1. Introduction}
\vspace{.2in}

A well-known theorem of Pringsheim says that for every
power series
\begin{equation}\label{1}
f(z)=\sum_{m=0}^\infty a_m z^m,\quad \limsup_{m\to\infty}|a_m|^{1/m}=1,
\end{equation}
with non-negative coefficients the point
$z=1$ is singular.

Fabry's theorem is a generalization of this;
assuming that projections of some coefficients on certain
lines through the origin
have relatively few sign changes, it guarantees
the existence of a singular point on a closed
arc of the unit circle
centered at $z=1$.

For the precise statement we need the following definitions.
For a sequence of real numbers $\{ a_m\}$, 
we say that a sign change occurs at the place $m$
if $a_ma_k<0$ for some $k<m$, while $a_{j}=0$ for $k<j<m$.

Let $\Lambda$ be a set of positive integers.
We denote by $n(r,\Lambda)$
the counting function
$$n(r,\Lambda)=\card\{\lambda\in\Lambda:\lambda\leq r\},$$
and define the {\em maximum density} of $\Lambda$ by the formula
$$D_2(\Lambda)=\lim_{r\to 0+}\limsup_{t\to\infty}
\frac{n((1+r)t,\Lambda)
-n(t,\Lambda)}{rt}.$$
The outside limit always exists, \cite[Satz III]{Polya1}. 
Here is an equivalent definition.
A set $\Lambda$ for which
$n(r,\Lambda)/r$ has a limit as $r\to\infty$ is called measurable,
and the limit is called the (ordinary) density of $\Lambda$.
Then $D_2(\Lambda)$ is the associated exterior density,
that is the infimum of densities of all
measurable sets of integers that contain $\Lambda$.
\vspace{.1in}

\noindent
{\bf Theorem A.} (Fabry \cite{Fabry})
{\em For a power series $f$ of the form $(1)$,
let $\{ m_k\}$ be a sequence with the property
\begin{equation}
\label{2}
\lim_{k\to\infty}|\Rea(e^{-i\beta_k}a_{m_k})|^{1/m_k}=1,
\end{equation}
with some real $\beta_k$.
Fix a number $r\in(0,1)$ and let $\Lambda_k$ be the set
of integers $m$ in the segment
\begin{equation}\label{3}
[(1-r)m_k,(1+r)m_k]
\end{equation}
where the sign chan\-ges of the se\-quen\-ce
$\{\Rea(e^{-i\beta_k}a_{m_k})\}$ occur.
If $\Delta=D_2(\cup_k\Lambda_k),$
then $f$ has a singularity on the arc
$$I_\Delta=\{ e^{i\theta}:|\theta|\leq\pi\Delta\}.$$}

The last sentence means that there is no immediate analytic
continuation of $f$ from the unit disc to the arc $I_\Delta$.
\vspace{.1in}

\noindent
{\em Comments.} 1. A sequence $\{ m_k\}$ satisfying (\ref{2})
always exists because the series $f$ has radius of
convergence $1$. One can take $\{ m_k\}$ such that
$|a_{m_k}|^{1/m_k}\to 1$ and then put $\beta_k=\arg a_{m_k}$.
Alternatively, one can first choose all
$\beta_k$ equal to $0$ or all $\beta_k$ equal to $\pi/2$,
and for at least one of these choices
a sequence $\{ m_k\}$ satisfying (\ref{2}) can be found.

2. Replacing $\{ m_k\}$ by a subsequence decreases\footnote{
Everywhere in this paper we use the words ``decrease'', ``increase''
etc. in the non-strict sense.} $\Delta$
and thus gives a stronger conclusion. For example, one can
add the assumption that intervals (\ref{3}) are disjoint,
and this will not weaken the result.

3. Same applies to the choice of the number $r$.
Choosing a smaller $r$ does not weaken the conclusion.
\vspace{.1in}

Fabry's statement in \cite{Fabry} is
equivalent to the statement above,
though he did not state a general
definition of the maximum density $D_2$.
This definition is due to P\'olya \cite{Polya1}. Bieberbach's
book \cite{Bieber}
contains a complete proof of Theorem A, 
as well as many corollaries and a survey of related
results up to the early 1950-s.
The history of Pringsheim's and
Fabry's
theorems
is described in \cite{Prings} by one of the main participants.
\vspace{.1in}

\noindent
{\bf Corollary 1.} {\em Suppose that for two real numbers
$\beta_j$, $0<\beta_1-\beta_2<\pi$
the set of sign changes in $\{\Rea(e^{-i\beta_j}a_n)\}$
has maximum density at most $\Delta$.
Then the power series $f$ in $(1)$ has a singularity on the 
arc $I_\Delta$.}
\vspace{.1in}   

The most often cited corollary of Fabry's theorem is this:
\vspace{.1in}

\noindent
{\bf Corollary 2.} {\em Suppose that the set $\{ n:a_n\neq 0\}$
has maximum density $\Delta$. Then $f$ has a
singularity on every closed arc of the unit circle of length
$2\pi\Delta$. In particular,
if $\Delta=0$ then the unit circle is the
natural boundary of $f$.}
\vspace{.1in}

Indeed, the number of sign changes of any sequence does not
exceed the number of its non-zero terms. So by Corollary 1
we conclude
that $f$ has a singularity on the arc $I_\Delta$.
But the conditions of Corollary 1 are invariant under a
transformation $f(z)\mapsto f(ze^{i\theta})$ so there is at
least one singularity on any closed
arc of length $2\pi\Delta$ on
the unit circle.

Other interesting corollaries are discussed
in the book \cite{Bieber}. Various special cases of
Theorem A were subject of intensive research in XX century,
however the fact that the assumptions of
the Theorem A can be substantially relaxed 
has been overlooked until recently.

One reason of this is that
Corollary 2 is best possible in a very strong sense
\cite[IX B]{Koosis}\footnote{Koosis credits
Fuchs \cite{Fuchs}
for the construction that proves this result.}:

\vspace{.1in}

{\em For every sequence $\Lambda$ of positive integers of maximum
density $\Delta>0$ and every $\delta\in (0,\Delta)$, there exists
a power series  $f$ of the form
$(\ref{1})$ with $a_n=0$ for $n\not\in \Lambda$, such that 
$f$ has an immediate analytic continuation from the unit disc
to the arc $\{ e^{i\theta}:|\theta|<\pi\delta\}.$}
\vspace{.1in}

In other words, the following two properties of a sequence
$\Lambda$ of positive integers are equivalent:
a) $D_2(\Lambda)\leq\Delta$ and b) every power series
of the form 
$$\sum_{m\in\Lambda}a_mz^m,\quad \limsup_{m\to\infty}|a_m|^{1/m}\to 1$$
has a singularity on the arc $I_\Delta$.
\vspace{.1in}

This result may create an impression that the 
maximal density is
the ``best possible density'' in Theorem A.
However we will see that this is not so. The difference
between Theorem A and corollaries 1, 2 is that
the density in Theorem A is measured not for the whole
sequence of coefficients but only for a part of it near
a subsequence $\{ a_{m_k}\}$ of ``large coefficients''. 

The first improvement of the density condition in Theorem A
is due to Arakelyan
and Martirosyan \cite{Arak}.
Suppose that a series (\ref{1}) and sequences $m_k,\beta_k$ 
satisfying (\ref{2}) are given. 
Let $\Lambda_{k,+}\subset[m_k,2m_k]$
and $\Lambda_{k,-}\subset[0,m_k]$ be the sets of
integers $j$ where the sign changes of $\Rea(e^{-\beta_k}a_j)$
occur. We denote $\Lambda_+=\{ \Lambda_{k,+}\}$ and
$\Lambda_-=\{\Lambda_{k,-}\}$, so that $\Lambda_{\pm}$
are sequences of finite sets of integers.
For every $r\in[0,1]$, we define 
\begin{equation}\label{maindef1}
n_{k,+}(r)=\frac{1}{m_k}\card\Lambda_{k,+}\cap [m_k,(1+r)m_k],
\end{equation}
and
\begin{equation}\label{maindef2}
n_{k,-}(r)=\frac{1}{m_k}\card\Lambda_{k,-}\cap [(1-r)m_k,m_k].
\end{equation}

Then we put
$$D_1(\Lambda_\pm)=\limsup_{r\to 0+}\limsup_{k\to\infty}
\frac{n_{k,\pm}(r)}{r}.$$

\noindent
{\bf Theorem B.} (Arakelyan and Martirosyan\footnote{
The statement given in \cite{Arak} is somewhat weaker,
but the argument there actually proves Theorem B.})
{\em The function $f$ 
in $(\ref{1})$ has a singularity on the arc $I_{\Delta}$,
where $\Delta=\min\{ D_1(\Lambda_+),D_1(\Lambda_-)\}$.}
\vspace{.1in}

According to Bieberbach,
P\'olya \cite[footnote 18 on p. 703]{Polya}
was the first to notice that in some versions of Fabry's theorem
the intervals (\ref{3}) can be replaced by
one-sided intervals $[m_k,(1+r)m_k]$ or by
$[(1-r)m_k,m_k]$. 

It is easy to see that
that $D_1(\Lambda_\pm)\leq D_2(\Lambda_\pm)$,
where $D_2(\Lambda_\pm)$
are the maximal densities of the sets
$\cup_k\Lambda_{k,\pm}$, 
and these inequalities can be strict. So Theorem~B
is stronger
than a ``one-sided'' version of Theorem A suggested by P\'olya.

The main result of this paper shows that the density condition
in Fabry's theorem can be further relaxed: we will replace $D_1$
by a smaller quantity.
To state it we need some preliminaries.
Notice that the functions $n_{k,\pm}$ defined in (\ref{maindef1})
and (\ref{maindef2}) are
increasing,
continuous from the right,
and satisfy the condition
\begin{equation}
\label{4}
|n_{k,\pm}(x)-n_{k,\pm}(y)|\leq|x-y|,\quad x,y\in[0,1],
\end{equation}
whenever $m_kx$ and $m_ky$ are integers.
By Helly's theorem, from every
sequence of such functions one can extract a
subsequence which
converges pointwise to some increasing function $n$.
We denote the sets of these limit functions $n$
by $\Fr(\Lambda_+)$ and $\Fr(\Lambda_-)$.
The limit functions
satisfy condition (\ref{4}) {\em for all real}
$x,y$ on $[0,1]$. In particular they are absolutely continuous
and their derivatives in the sense of distributions
satisfy $\| n'\|_\infty\leq 1.$
We also have $n(0)=0$.   

Let $n$ be an increasing function on some closed interval $I$
of the real line,
satisfying the condition $|n(x)-n(y)|\leq|x-y|$ for all $x,y$ in $I$.
For every $\Delta\in[0,1]$ we define
the {\em lower $\Delta$-regularization},
$$\underline{n}^\Delta=\underline{n}^\Delta_I
=\sup\{\phi\in C^1:\phi\leq n,\;
\Delta\leq\phi'\leq 1\}.$$
In other words, $\underline{n}^\Delta_I$
is the largest minorant of $n$ on $I$
whose slope is at least $\Delta$. Notice that
$\Delta_1<\Delta_2$ implies $\underline{n}^{\Delta_1}_I\geq
\underline{n}^{\Delta_2}_I$,
and $\underline{n}^0_I=n$. Furthermore, if $n$ is originally
defined on $I$ and we restrict 
it to a smaller interval $I_1\subset I$,
and take a regularization of this restriction, then
$\underline{n}^\Delta_{I_1}(x)\geq \underline{n}^\Delta_{I}(x)$
for $x\in I_1$.

If $I=[0,\delta], \delta>0$, 
we will denote $\underline{n}^\Delta_{I}$ by
$\underline{n}^\Delta_{\delta}$. In what follows we will
sometimes simplify the notation by omitting any reference to
the interval of regularization, if this interval
is clear from context. When doing this, we will always use the
following convention: in a regularization that occurs
in an integrand,
the interval of the regularization
coincides with the interval
of integration.
\vspace{.1in}

{\bf Theorem 1.} {\em Let a power series $f$ as in $(\ref{1})$
and sequences
$\{ m_k\}, \{\beta_k\}$ with the property $(\ref{2})$ be given.
If for some limit function
 $n\in\Fr(\Lambda_+)\cup\Fr(\Lambda_-)$
and a number
$\Delta\in[0,1)$ we have
\begin{equation}\label{5}
\int_{0}^\delta\frac{n(r)-\underline{n}^\Delta_\delta(r)}{r^2}dr
=\infty, \quad\mbox{for all}\quad\delta\in(0,1),
\end{equation}
then $f$ has a singularity on the arc $I_{\Delta}$.}
\vspace{.1in}

This result suggests the following definitions:
$$D_3(\Lambda_+)=\inf_{n\in\Fr(\Lambda_+)}
\inf\left\{ a\in[0,1]:\;\int_0^\delta
\frac{n(r)-\underline{n}^a_\delta(r)}{r^2}dr=\infty,\;
\mbox{for all}\;\delta\in(0,1)\right\},$$
and similarly for $\Lambda_-$. 

We will show (Lemma 3 in section 2) that the densities
$D_3$ have the following monotonicity property:
if for every $k$ we have $\Lambda_{k,\pm}\subset
\Lambda^\prime_{k,\pm}$,
then $D_3(\Lambda_\pm)\leq D_3(\Lambda^\prime_{\pm})$.

This monotonicity property combined with Theorem 1
gives a ``gap version'' of Theorem 1:
instead of counting sign changes we can define
$\Lambda_\pm$
as the sequences of subscripts $j$ of {\em non-zero terms}
of 
$\{ a_j\}$ for $m_k\leq j\leq 2m_k$ and
$0\leq j\leq m_k$ respectively. Then
$f$ has at least one singularity on each closed arc
of the unit circle of length $\pi\Delta$, 
where $\Delta=\min\{ D_3(\Lambda_+),D_3(\Lambda_-)\}$.

To compare Theorem 1 with theorems A and B, we choose
the sequence $\{ m_k\}$ in (\ref{2}) in
such a way that the limit
$n=\lim_{k\to\infty}n_{k,+}$ exists.
Replacing $\{ m_k\}$ by its subsequence
can only decrease the densities $D_2$ and $D_1$.
Then 
$$D_1(\Lambda_+)\geq\limsup_{r\to 0+}n(r)/r,\quad\mbox{and}\quad
D_2(\Lambda_+)\geq \limsup_{r,r'\to 0+}|n(r)-n(r')|/|r-r'|.$$
On the other hand, it is easy to see that
$$D_3(\Lambda_+)\leq\liminf_{r\to 0+}n(r)/r.$$
Let us combine this with Theorem 1 to obtain a corollary
whose conditions
are easier to verify:
\vspace{.1in}

\noindent
{\bf Corollary 3.} {\em Let a power series $f$ as in $(\ref{1})$ and
sequences $\{ m_k\}, \{\beta_k\}$ satisfying $(\ref{2})$ be given.
If some limit function $n\in\Fr(\Lambda_+)\cup\Fr(\Lambda_-)$
satisfies
$$\liminf_{r\to 0+}\frac{n(r)}{r}\leq\Delta,$$
then  $f$ has a singularity on the arc $I_\Delta$.}
\vspace{.1in}

We summarize the relations between the considered densities as 
$$ D_3\leq D_1\leq D_2,$$
and all inequalities can be strict. Assuming that
$\Fr(\Lambda_+)\cup\Fr(\Lambda_-)$ contains a function $n$
such that 
$$\limsup_{r\to 0+}n(r)/r=1\quad\mbox{and}\quad
\liminf_{r\to0+}n(r)/r=0,$$
we obtain $D_3=0$ while $D_1=D_2=1$.
In this case, Theorems A and B say nothing, while Theorem 1
implies that $z=1$ is a singular point, and the gap version
of Theorem 1 
gives that the
whole unit circle is the natural boundary.  

In the recent paper \cite{Ar}, a new
density condition in Fabry's gap theorem is given, which is
incomparable with our conditions in
Theorem 1 or its Corollary 3.
The density used in \cite{Ar}
can be written in our notation
as
$$D_4(\Lambda)=\liminf_{r\to 0}\liminf_{k\to\infty}
\frac{1}{2r}\int_0^r\frac{n_{k,+}(t)+n_{k,-}(t)}{t}dt,$$
and it is shown that every power series (\ref{1})
with
$$|a_{m_k}|^{1/m_k}\to 1$$
and $a_j=0$ for $j\in[0,2m_k]\backslash\Lambda_k$,
has a singularity
on the arc $I_{D_4(\Lambda)}$.

If the pointwise limits
$n_{\pm}=\lim_{k\to\infty}n_{k,\pm}$
exist, then it is easy
to see that $$D_4(\Lambda)\geq
\min\{ \liminf_{r\to 0}n_+(r)/r,\liminf_{r\to0}
n_-(r)/r\}.$$
So in this case, our Corollary 3 
gives a stronger result.

However one can construct examples in which
$$D_4(\Lambda)<\min\{ D_3(\Lambda_+),D_3(\Lambda_-)\},$$
so in general our Theorem 1 does not contain the
result of \cite{Ar} as a special case.

Thus the question on the best possible density condition
in Fabry's theorem remains open.  
\vspace{.1in}

{\em Sketch of the proof of Theorem 1.}

Assume for simplicity that the coefficients $a_m$ are real
and choose $\beta_k=0$. If $f$ has an immediate analytic
continuation on $I_\Delta$, then the sequence $(-1)^ma_m$ can be
interpolated by a holomorphic function $F$ in some angle containing
the positive ray, such that
$\log|F(z)|\leq \pi b|\Ima z|+o(|z|),\; z\to\infty$
(Theorem C in section 2). If the sequence $\{ a_m\}$
has few sign changes
on some interval, then $F(m)=(-1)^ma_m$ has many zeros on the
same interval
(Lemma~1, section 2).
Thus we need to estimate from above the number of zeros of $F$
near the points $m_{k}$ where $|F(m_k)|=|a_{m_k}|$ is not too small
(is it not too small by (\ref{2})).
After a more or less standard rescaling trick, this is reduced
to an estimate from above of the Riesz measure of a subharmonic
function $u$ in a neighborhood of $0$ having the
properties $u(0)=0$ and $u(z)\leq \pi b|\Ima z|$.
Such estimate of a Riesz measure from above can be obtained
by adapting the arguments of Beurling and Malliavin from
\cite{BM} (Lemmas 2--6 in section 2).

Now we give the details.
\vspace{.2in}

\newpage
\noindent
{\bf Preliminary results.}
\vspace{.2in}

\nopagebreak
We will use the following
\vspace{.1in}

\noindent
{\bf Theorem C.} {\em For a function $f$ as in $(\ref{1})$ to
have an immediate analytic continuation from the unit disc
to the arc
$I_\Delta$ it is necessary and sufficient
that there exists a function $F$ analytic in some angle
$A(\alpha)=\{ z:|\arg z|<\alpha\}$ with the properties
\begin{equation}
\label{7a}
a_m=(-1)^mF(m),
\end{equation}
and 
\begin{equation}
\label{7b}
\limsup_{t\to\infty}\frac{\log|F(te^{i\theta})|}{t}
\leq \pi b|\sin\theta|,\quad |\theta|<\alpha,
\end{equation}
with some $b<1-\Delta.$}
\vspace{.1in}

This is a special case of
\cite[Ch. V, Th. III]{Bern}\footnote{There is a misprint
in Bernstein's statement: his inequality (14) should
be $|t|\leq\ell$. With $|t|<\ell$, Bernstein's
statement no longer holds,
even for power series. See
\cite{Ar}.}.
A simple proof of this special case can be found in
\cite{Arak}. We only need the ``necessary'' part of this
theorem, and we include a proof for the reader's convenience.
\vspace{.1in}

{\em Proof of necessity.} We begin with a function
$F_\epsilon$ defined by the formula
$$F_\epsilon(z)=\frac{1}{2\pi i}\int_{-i\pi-\epsilon}^{i\pi-\epsilon}
f(-e^{\zeta})e^{-z\zeta}d\zeta,$$
where 
$\epsilon>0$ is arbitrary. Then $F_\epsilon$ is an entire function of exponential
type.
Cauchy's formula gives
$$(-1)^ma_m=\frac{1}{2\pi i}\int_{|w|=\rho} f(-w)w^{-m-1}dw,
$$
where $\rho<1$.
Making the change of the variable $w=e^{\zeta}$ in the Cauchy
integral, we obtain (\ref{7a}) for all functions
$F_\epsilon$. 

The integrand in $F_\epsilon$
is analytic in the left half-plane,
and by assumption it has an immediate analytic continuation
to a neighborhood of the two segments $[-i\pi,-i\pi b]$
and $[i\pi b,i\pi]$
for some $b<1-\Delta$.
So we can deform the path of integration to a new
path $\gamma$ shown in Fig. 1.
\vspace{.1in}
\begin{center}
\epsfxsize=3.0in
\centerline{\epsffile{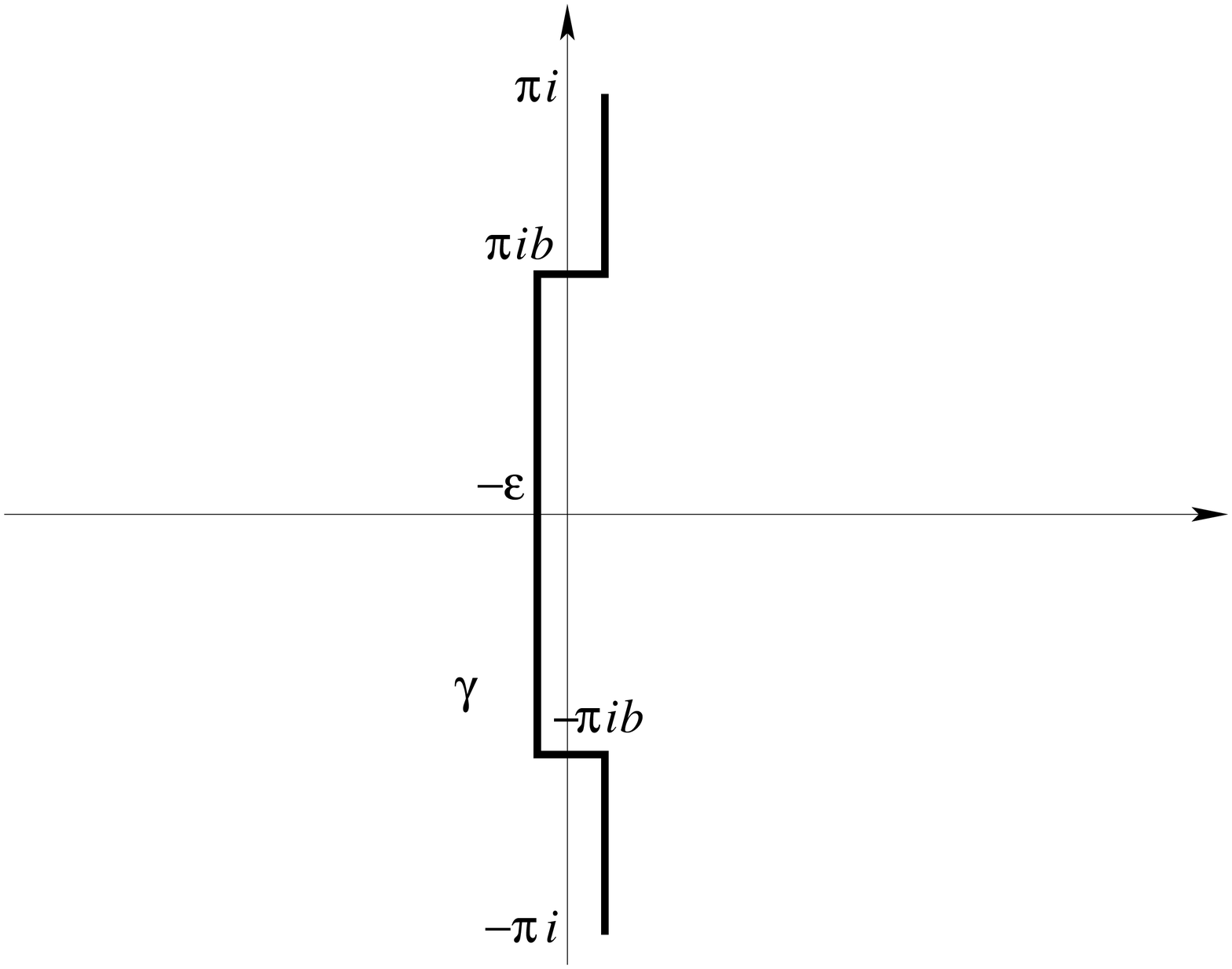}}
\nopagebreak
\vspace{.1in}
Fig. 1. Path $\gamma$.
\vspace{.1in}
\end{center}
This path $\gamma$ consists of the vertical segment
$[-i\pi b-\epsilon,i\pi b-\epsilon]$, two horizontal segments
$[\pm i\pi b-\epsilon,\pm i\pi b+\epsilon_1]$ and two
vertical segments on the line $\Rea\zeta=\epsilon_1$.
This path deformation changes $F_\epsilon$,
but does not change its values at the positive integers,
because when $z$ is an integer, the integrand in $F_\epsilon$
has period $2\pi i$. 
Now we set
$$F(z)=\frac{1}{2\pi i}\int_\gamma f(-e^\zeta)e^{-z\zeta}d\zeta,$$
and $F$ satisfies (\ref{7a}).
The function 
$$h_\gamma(\theta)=\sup_{\zeta\in\gamma}(-\Rea(\zeta
 e^{i\theta}))$$
satisfies $$h(\theta)=\pi b|\sin\theta|+\epsilon\cos\theta,
\quad |\theta|\leq\alpha,$$
for some $\alpha>0$ depending only on $\epsilon_1$.
Then the straightforward
estimate of the integral
over $\gamma$ gives
$$|F(z)|\leq C \exp(|z| h(\arg z)),\quad |z|\to\infty$$
where $C$ is a constant depending on $f$ (see, for example,
\cite[Ch. I, \S 19-20]{Levin}).
Changing $\epsilon$ does not change $F$
by Cauchy's theorem.  
 Letting $\epsilon\to 0+$,
we obtain  
(\ref{7b}).
\hfill$\Box$
\vspace{.1in}

\noindent
{\bf Lemma 1.} {\em Let $(a_0,a_1,\ldots,a_N)$ be a sequence of
real numbers, and $f$ a real analytic function on the closed interval
$[0,N]$, such that $f(n)=(-1)^na_n$. Then the number of zeros
of $f$ on $[0,N]$, counting multiplicities, is at least
$N$ minus the number of sign changes of the sequence $\{ a_n\}$.}
\vspace{.1in}

{\em Proof.} Consider first an interval $(k,n)$ such that
$a_ka_n\neq 0$ but $a_j=0$ for $k<j<n$.
We claim that $f$ has at least 
$$n-k-\#(\mbox{sign changes in the pair}\; (a_k,a_n))$$
zeros
on the open interval $(k,n)$. Indeed, the number of zeros of $f$
on this interval is at least $n-k-1$ in any case.
This proves the claim if there is a sign change in the pair $(a_k,a_n)$. 
If there is no sign change, that is $a_na_k>0$, then  $f(n)f(k)=(-1)^{n-k}$.
So the number of zeros of $f$ on the interval $(n,k)$ is of the same
parity as $n-k$. But $f$ has at least $n-k-1$ zeros on this
interval, thus the total number of zeros is at least $n-k$.
This proves our claim.

Now let $a_k$ be the first and $a_n$ the last non-zero term of our sequence.
As the interval $(k,n)$ is a disjoint union of the intervals
to which the above claim applies, we conclude that the number of
zeros of $f$ on $(k,n)$ is at least $(n-k)$ minus the number of
sign changes of our sequence. On the rest of the interval 
$[0,N]$
our function has at least $N-n+k$ zeros, so the total
number of zeros is at least $N$ minus the number of sign changes.
\hfill$\Box$
\vspace{.2in}

Let $n$ be an increasing function on a closed interval $I$.
For every $a\in [0,1]$ we define the {\em upper 
$a$-regularization} by
$$\overline{n}^a_I=\inf\{ \phi:\phi\geq n,\;
 0\leq\phi'\leq a\}.$$ 
For functions $n$ satisfying (\ref{4}) we have the formula
\begin{equation}
\label{formula}
\overline{(\id-n)}^a_I=\id-\underline{n}^{1-a}_I,
\end{equation}
which is easy to verify.
\vspace{.1in}

Consider the set $\N$ of all increasing functions $n, n(0)=0$
on a segment $I=[0,\delta]$,
where $\delta>0$ is fixed. We introduce the following order relation
$n_1\succ n_2$ if $n_1-n_2$ is increasing.
\vspace{.1in}

\noindent
{\bf Lemma 2.} {\em For $n\in\N$: 
\begin{equation}\label{A}
\int_0^\delta\frac{\overline{n}^a(r)-n(r)}{r^2}dr<\infty
\end{equation}
if and only if there exists $n_1\in\N$ with the properties
$n_1\succ n,\; n_1(r)\leq ar,\; 0\leq r\leq\delta$ and
\begin{equation}\label{B}
\int_0^\delta\frac{ar-n_1(r)}{r^2}dr<\infty.
\end{equation}}

{\em Proof.} (\ref{A})$\longrightarrow$(\ref{B}).
Put $n_1(r)=ar-\overline{n}^a+n.$ It satisfies all conditions.
\vspace{.1in}

\noindent
(\ref{B})$\longrightarrow$(\ref{A}). We define $n_2(r)=n(r)+ar-n_1(r),$
then $n_2\geq n$ and $a\cdot\id \succ n_2$. This implies that
$n\leq\overline{n}^a\leq n_2$, and by (\ref{B})
$$\int_0^\delta \frac{n_2(r)-n(r)}{r^2}dr<\infty,$$
holds. We conclude that (\ref{A}) holds as well.
\hfill$\Box$
\vspace{.1in}

{\bf Lemma 3.} {\em If $n$ and $n_1$ are in $\N$, and $n_1\succ n$, then
\begin{equation}\label{AA}
\int_0^\delta\frac{\overline{n}^a(r)-n(r)}{r^2}dr=\infty
\end{equation}
implies
\begin{equation}\label{BB}
\int_0^\delta\frac{\overline{n_1}^a(r)-n_1(r)}{r^2}dr=\infty.
\end{equation}
}

{\em Proof.} Suppose that the integral in (\ref{BB}) converges. 
By Lemma 2 there exists $n_2\succ n_1,\; n_2\leq a\cdot\id$
such that 
$$\int_0^\delta\frac{ar-n_2(r)}{r^2}dr<\infty.$$
As $n_2\succ n_1\succ n$ and $n_2\leq a\cdot\id$, another application
of Lemma 2 yields that the integral in (\ref{AA}) converges.
\hfill$\Box$
\vspace{.1in}

\noindent
{\bf Lemma 4.} {\em Let $u$ be a subharmonic function
in $\{ z:|z|<2\delta\}$, satisfying
\begin{equation}\label{prop1}
u(0)=0,
\end{equation}
and
\begin{equation}\label{prop2}
u(z)\leq \pi b|\Ima z|,\quad |z|<2\delta,
\end{equation}
for some $b>0$.
Then
\begin{equation}
\label{prop3}
\int_{-\delta}^{\delta}\frac{u(x)}{x^2}dx>-\infty.
\end{equation}}

{\em Proof.} We may assume without loss of generality
that $u(z)=u(\overline{z})$ (replacing $u$ by $(u(z)+u(\overline{z}))/2$
alters neither the conditions nor the assumptions of the lemma).
Consider the Poisson integral in the upper half-plane 
$$v(x+iy)=\frac{y}{\pi}\int^\delta_{-\delta}
\frac{u(t)}{(x-t)^2+y^2}dt.$$
This integral is convergent because $u$ is intergable
on the interval $(-\delta,\delta)$.
Let $w$ be the least harmonic majorant for the
subharmonic function $u-v$ in the half-disc
$D=\{ z:|z|<\delta,\Ima z>0\}$.
Then $w$ is a harmonic function in $D$,
whose limit on the diameter of $D$ is zero.
By reflection, $w$ extends to a harmonic function
in the whole disc $\{ z:|z|<\delta\}$.
It follows that the normal derivative $\partial w/\partial y$
is bounded on the interval $-\delta/2<x<\delta/2$.
So there exists a neighborhood $V$ of $0$ and a constant
$c>0$ such that 
$$u(z)\leq v(z)+\pi c|\Ima z|,\quad z\in V,\quad\Ima z>0.$$
Suppose that the integral in (\ref{prop3}) is divergent,
then $v(iy)/y\to-\infty$ and thus
$u(iy)/y\to-\infty$ as $y\to 0+$. 
Thus there exists $y_0>0$ such that 
\begin{equation}\label{prop4}
u(iy)\leq -y,\quad 0\leq y\leq y_0.
\end{equation}
Now we consider the sequence of subharmonic functions
$u_n(z)=2^nu(2^{-n}z).$ By (\ref{prop2}), this sequence is
uniformly bounded from above on compact subsets of the plane,
and by (\ref{prop1}) it is bounded from below at $0$.
Compactness Principle \cite[Th. 4.1.9]{Hor}
implies that some subsequence of $\{ u_n\}$
converges in $L^1_{\rm loc}$
to a function $u_{\infty}$ subharmonic in the whole plane.
Moreover, 
$$\limsup_{n\to\infty}u_n(z)\leq u_\infty(z),\quad z\in\C,$$
by the same theorem in \cite{Hor}.
In view of (\ref{prop2}), this function $u_\infty$ satisfies 
$u_\infty(z)\leq \pi b|\Ima z|$ in the whole plane, and 
in addition it follows from (\ref{prop4})
$u_\infty(iy)\leq -|y|$ for all real $y$. Here we used
the symmetry assumption made in the beginning of the proof.
These two properties contradict the Phragm\'en--Lindel\"of Principle,
which proves the lemma.\hfill$\Box$
\vspace{.1in}

\noindent
{\bf Lemma 5.} {\em Let a countable set of open intervals
whose lengths tend to zero be given, and let $E$ be the union
of these intervals.
Then there exists a subset of these intervals whose union
is also $E$, and no point of $E$ belongs to more than
two intervals of the subset.}
\vspace{.1in}

{\em Proof.} We order the given intervals into a sequence of
decreasing length. Inspecting the intervals of this sequence
one after another, we select or discard them.
On the first step, the first interval is selected. On the $k$-th step,
the $k$-th interval of the sequence is discarded
if it belongs to the union of the intervals selected on
the previous steps, otherwise this $k$-th interval is selected.

Consider now all selected intervals. It is clear that their
union is $E$, because on every step the
union of non-discarded intervals does not change.

We claim that every point of $E$ is covered by
finitely many selected intervals. Indeed, let $x$ be a point of $E$.
Let $I$ be some selected interval containing $x$. Suppose that
$I$ was selected on $k$-th step. 
If $x$ is covered by infinitely many selected intervals,
infinitely many of them are contained in $I$
because the lengths of the
intervals tend to zero. Then some of these infinitely many intervals
containing $x$ had to be selected after step $k$, which contradicts
the selection rule. This proves the claim.

Now we remove all those selected intervals
which are contained in the union
of other selected intervals. We claim that the intervals
that were not removed still cover $E$. Indeed, let $x$ be a point
in $E$.
Then $x$ belongs to finitely many selected intervals. And it is
evidently impossible that each interval of a finite family of intervals
is contained in the union of the rest.  

So the remaining intervals have the property that none of them
is contained in the union of the rest.
Such family of intervals cannot have
triple intersections: if three intervals intersect, then one of them
is contained in the union of the other two.
\hfill$\Box$
\vspace{.1in}

In the following lemma we will have to deal with restrictions
of increasing functions $\nu$ to smaller intervals.
We recall that if  we restrict $\nu$ to a smaller interval $I'\subset I$,
the upper $a$-regularization of this restriction
will be less than or equal to the restriction to $I'$ of the 
upper $a$-regularization of  $\nu$ on $I$.
If $I=[0,\eta]$ we write $\overline{\nu}^a_\eta$
instead of $\overline{\nu}^a_I$.
\vspace{.1in}

\noindent
{\bf Lemma 6.} {\em Let $u$ be a function from Lemma 4.
Denote by $\nu(r)$ the Riesz measure corresponding to $u$
of the segment $[0,r]$.
Then for every $a>b$ there exists $\eta\in (0,\delta)$ such that
\begin{equation}
\label{main}
\int_{0}^\eta\frac{\overline{\nu}^a_\eta(r)-\nu(r)}{r^2}dr<\infty.
\end{equation}}

{\em Proof.} We follow Kahane's exposition \cite{Kahane} of the
work of Beurling--Malliavin \cite{BM}.
Jensen's formula and (\ref{prop2}) give
\begin{eqnarray*}
u(x)&\leq&-\int_0^R (\nu(x+t)-\nu(x-t))\frac{dt}{t}+
\frac{1}{2\pi}\int_{-\pi}^\pi u(x+Re^{i\theta})d\theta\\
&\leq&-\int_0^R (\nu(x+t)-\nu(x-t))\frac{dt}{t}+2bR.
\end{eqnarray*}
Integrating this with respect to $x$ from $\alpha$ to $\beta$,
and using the estimates
$$\int_{\alpha-t}^{\alpha+t}\nu(x)dx\leq 2t\nu(\alpha+R)$$
and
$$\int_{\beta-t}^{\beta+t}\nu(x)dx\geq 2t\nu(\beta-R),$$
which follow from monotonicity of $\nu$,
we obtain
\begin{equation}\label{6}
\int_\alpha^\beta u(x)dx\leq2R(b(\beta-\alpha)-(\nu(\beta-R)-\nu(\alpha+R))).
\end{equation}
Suppose now that for some interval $(\alpha,\alpha+\ell)$ we have
$\nu(\alpha+\ell)-\nu(\alpha)\geq a\ell.$
Putting 
\begin{equation}\label{epsilon}
\epsilon=(b-a)/(2(b+a)),\quad \beta=\alpha+\ell+\ell\epsilon,\quad R=\ell
\epsilon,
\end{equation}
we obtain from (\ref{6}) that
\begin{equation}
\label{7}
\int_\alpha^\beta u(x)dx\leq-2\epsilon\epsilon'\ell^2,
\end{equation}
where $\epsilon'=a(1-\epsilon)-b(1+\epsilon)>0.$


The set $E=\{ x:\overline{\nu}^a_\delta(x)>\nu(x)\}$ consists
of disjoint open intervals $J_n=(\alpha_n,\alpha_n+\ell_n).$
We may assume that the union of these
intervals has $0$ as an
accumulation point, otherwise (\ref{main}) holds trivially.
\vspace{.1in}

Case 1. Suppose that $0$ is not an endpoint of any interval $J_n$.
Then
\begin{equation}
\label{aa}
\int_0^\delta
\frac{\overline{\nu}^a_\delta(x)-\nu(x)}{x^2}dx=
\sum_n\int_{J_n}
\frac{\overline{\nu}^a_\delta(x)-\nu(x)}{x^2}dx\leq
\sum_n\frac{\ell_n^2}{\alpha_n^2}.
\end{equation}
The enlarged intervals $J_n^\prime=(\alpha_n,\beta_n),$
where $\beta_n=\alpha_n+\ell_n+\epsilon\ell_n$
might no longer be disjoint,
but we can apply Lemma 5 to find a subset of
these intervals that
covers $E$ with multiplicity at most $2$.
Then, using (\ref{7}), we obtain
\begin{equation}
\label{17}
-\infty<\int_0^\delta \frac{u(x)}{x^2}\leq
2\sum_n\frac{1}{\beta_n^2}\int_{\alpha_n}^{\beta_n} u(x)dx\leq-
 4\epsilon\epsilon'\sum_n
\frac{\ell_n^2}{\beta_n^2},
\end{equation}
so the last series converges.
But then $\ell_n/\beta_n\to 0$, so $\alpha_n\sim\beta_n$,
and we conclude that the series in the right hand side of
(\ref{aa}) also converges. This proves the lemma with $\eta=\delta$
in this case.
\vspace{.1in}

Case 2. Suppose now that some interval $J$ has the form $J=(0,x_0)$.
Then $\nu(x_0) \geq ax_0$. We may decrease the interval $[0,\delta]$
on which the majorant is defined, and perhaps obtain a new majorant
$\overline{\nu}^a_\eta$ on a smaller interval $[0,\eta]$, such that
the new set $E=\{ x:\overline{\nu}^a_\eta(x)>\nu(x)\}$
will not contain an interval $J$ with an endpoint at $0$.
Then we repeat the argument of the Case 1.

Otherwise, there is a sequence $x_k\to 0$ such that $\nu(x_k)\geq ax_k$,
and the majorants $\overline{\nu}^a_{x_k}$ on $[0,x_k]$ have
the property $\overline{\nu}^a_{x_k}(x)>\nu(x)$ for $x\in(0,x_k).$
In particular, $\nu(x_k)-\nu(x_k/2)\geq ax_k/2.$
We can choose a subsequence so that the intervals
$(x_k/2,2x_k)$ are disjoint.
Taking $\alpha_k=x_k/2$, $\ell_k=x_k/2$, and
$\beta_k=x_k+\epsilon x_k/2$, where $\epsilon$ is defined in (\ref{epsilon}),
we obtain
intervals to which the inequality (\ref{7}) applies,
so we can write (\ref{17}) again, and obtain a contradiction
because this time $\ell_n/\beta_n$ does not tend to zero.
\hfill$\Box$

%
%
\vspace{.2in}

\noindent
{\bf Proof of Theorem 1.}
\vspace{.2in}

Proving the theorem by contradiction,
we will assume that (\ref{5}) holds for a limit function
$n$ of $n_{k,+}$, and that
$f$ has an immediate analytic
continuation through the arc $I_\Delta.$ 
The case of a limit function of $n_{k,-}$
is completely similar.

Applying Theorem C to $f$ we obtain a function $F$ holomorphic
in some angle $A(\alpha)$ with
the properties (\ref{7a}) and (\ref{7b}). Assume that for our
sequence $\{ m_k\}$ the limit $\lim_{k\to\infty} n_k=n$
satisfying (\ref{5}) exists.
Consider the sequence 
$$F_k(z)=e^{-i\beta_k}F(z)+e^{i\beta_k}\overline{F(\overline{z})}.$$
These functions are real on the positive ray, and satisfy
\begin{equation}
\label{18}
F_k(m)=2(-1)^m\Rea(a_me^{-i\beta_k}),
\end{equation}
thus by Lemma 1, the number of zeros of $F$ on
every interval $(m',m^{\prime\prime})\subset[m_k,2m_k]$
with integer endpoints is at least
\begin{equation}\label{MMM}
m^{\prime\prime}-m^\prime-\#(\mbox{changes of sign $\{ \Rea(a_je^{i\beta_k})\}$
for 
$m'\leq j\leq m"$}).
\end{equation}
Consider the subharmonic functions
$$u_k(z)=\frac{1}{m_k}\log|F_k(m_k(z+1))|.$$
In view of (\ref{7b}) this sequence of subharmonic functions
is uniformly bounded from above on every compact subset of
the angle $A(\alpha)-1$. Moreover, condition (\ref{2}) together
with (\ref{18}) imply that the $u_k(0)$ are bounded from below.
Then the Compactness Principle for subharmonic functions
\cite[Th. 4.1.9]{Hor} implies that, after choosing a subsequence,
 $u_k\to u$, where $u$ is a subharmonic function
in the angle $A(\alpha)-1$. This function $u$ has the properties
(\ref{prop1}) and (\ref{prop2}) of Lemma 4 with $b<1-\Delta$,
if $\delta<\sin\alpha$. 
Choose $a\in(b,1-\Delta)$.
The Riesz measures of $u_k$ converge to the Riesz
measure of $u$ weakly. 
Let $\nu(r)$ be the Riesz measure corresponding to $u$
of the interval $[0,r]$. Then (\ref{MMM}) implies that $\nu\succ \id-n$.
Using Lemma 3 and (\ref{formula}), we conclude
that that for every $\eta\in (0,\delta)$
$$\int_0^\eta\frac{\overline{(\id-n)}^a_\eta(r)-r+n(r)}{r^2}dr=\infty.$$
Now Lemma 3 implies that
$$\int_0^\eta
\frac{\overline{\nu}^a_\eta(r)-\nu(r)}{r^2}dr=\infty,$$
and this contradicts Lemma 6. 
\hfill$\Box$
\vspace{.1in}

I thank Alan Sokal for many interesting conversations
that stimulated this research,
Andrei Gabrielov for the proof of Lemma 5 and
useful comments and Norair Arakelian who brought
\cite{Ar} to my attention and the referee whose
suggestions improved the exposition.

\vspace{.1in}

{\em Purdue University

West Lafayette, IN 47907-2067 USA

eremenko@math.purdue.edu}

\end{document}